\documentclass[12pt]{article}

\usepackage{amsmath}
\usepackage{amssymb}

\newcommand{\eop}{\bigstar}  

\newenvironment{Proof}{\noindent{\bf Proof.}}{\par\bigskip} 

\newenvironment{Proof of the Claim}{\noindent{\bf Proof of the Claim.}}{\par\bigskip}

\newtheorem{Theorem}{Theorem}[section]
\newtheorem{Lemma}[Theorem]{Lemma}
\newtheorem{Definition}[Theorem]{Definition}

\newtheorem{Claim}[Theorem]{Claim}

\newcommand{\elementary}{\prec}

\newcommand{\deq}{\buildrel{\rm def}\over =}

\newcommand{\FF}{{\cal F}}

\newcommand{\KK}{{\cal K}}

\newcommand{\TT}{{\cal T}}

\newcount\skewfactor
\def\mathunderaccent#1#2 {\let\theaccent#1\skewfactor#2
\mathpalette\putaccentunder}
\def\putaccentunder#1#2{\oalign{$#1#2$\crcr\hidewidth
\vbox to.2ex{\hbox{$#1\skew\skewfactor\theaccent{}$}\vss}\hidewidth}}

\begin{document}

\title{A general Stone representation theorem}

\author{Mirna; after a paper by A. Jung and P. S\"underhauf and\\ notes by G. Plebanek}

\maketitle

This note contains a Stone-style representation theorem for compact Hausdorff spaces.

The note is very much inspired by some existing representation theorems and is expository
in nature. The
first representation theorem is by
Jung and S\"underhauf in \cite{JS} and there is also a version of it for compact Hausdorff
spaces noted by D. Moshier \cite{Mo}, the ideas of which were
described to us by A. Jung. This covers \S{\ref{spils} of the note. The other
representation theorem uses normal lattices and was discovered recently by G. Plebanek
\cite{Grz}. We show in \S2 that the two representation theorems are reducible to
each other.

Some notation and definitions might be specific to this note.
The original, highly recommended paper \cite{JS}, deals with a much more general
situation where neither the Hausdorff property nor compactness
is assumed in the representation theorem, while the notes \cite{Grz}
develop the lattice representation theorem in a self-sufficient manner.
Further results on representation theorems
likely including the theorem described in \S1 here, will be published as part of a future
paper by Jung et al, and the work of Plebanek is done as part of a separate project on Banach
spaces.

\section{Spils}\label{spils}

\begin{Definition}\label{sba} A {\em strong proximity
involution lattice (spil)} is given by a structure
$\langle B, \vee, \wedge,\,{}^{'}, 0,1,\elementary\rangle$ where 
$\langle B, \vee, \wedge, 0,1\rangle$ is a bounded distributive lattice
and the following
additional axioms hold:
\begin{description}
\item{(i)}\label{i}
$\elementary$ is a binary relation which is {\em interpolating},
meaning it satisfies $\prec^2=\prec$ so for all $a,b,c\in B$
\begin{description}
\item{(a)} $a\prec b\,\,\&\,\,b\prec c\implies a\prec c$
\item{(b)} if $a\prec c$ then there is some $b$ such that $a\prec b
\,\,\&\,\, b\prec c$;
\end{description}
\item{(ii)}\label{ii}
\begin{description}
\item{(a)} for all finite $M\subseteq B$ and $a\in B$
\[
M\elementary a\iff \bigvee M\elementary a;
\]
\item{(b)}\label{iii} for all finite $M\subseteq B$ and $a\in M$
\[
a\elementary M\iff a\elementary  \bigwedge M;
\]
\end{description}
\item{(iii)}\label{vi} Involution ' is a unary operation satisfying that
\begin{description}
\item{(a)} $x''=x$ for all $x$ (so the involution is proper);
\item{(b)}
for all $x,y$ and $z$ we have $x \wedge y \prec z$  iff  $x \prec z \vee y'$
and
\item{(c)} (De Morgan laws) $(x\vee y)'=x'\wedge y'$ and its dual $(x\wedge y)'=x'\vee y'$ hold; 
\end{description}
\item{(iv)} $x \prec  y \wedge y'  \implies x \prec 0$.
\end{description}
\end{Definition}

Here we use the notation $M\prec a$ for $(\forall m\in M)\, m\prec a$ and 
similarly for
$a\prec M$.
The idea of a spil is that it is a substitute for a Boolean algebra,
where the involution plays the role of the complement and 
$\prec$ the role of the order $\le$ induced by the Boolean operations.
As in the classical case of the Boolean algebras there is 
some duality in the axioms, as indicated by pairing (a) and (b) in (ii) and (iii).

Some basic properties of spils are given by the following Lemma, which
is Lemma 7 in \cite{JS}. For the sake of completeness we give the proof
here as well.

\begin{Lemma}\label{basic}
Suppose that $B$ is a spil. Then for all $a,b, c,d\in B$ we have
\begin{description}
\item{(1)} $0\prec a\prec 1$,
\item{(2)} $a\prec b\implies a\prec b\vee c$,
\item{(3)} $a\prec b\implies a\wedge c\prec b$,
\item{(4)} $a\prec b\,\,\&\,\,c\prec d\implies a\vee b \prec c\vee d$,
\item{(5)} $a\prec b\,\,\&\,\,c\prec d\implies a\wedge b\prec c\wedge d$.
\end{description}
\end{Lemma}

\begin{Proof} (1) We have $\emptyset\prec a$ trivially so $0=\bigvee\emptyset
\prec a$ by axiom (ii)(a) of a spil. Similarly $a\prec \bigwedge \emptyset=1$ by
(ii)(b). For (2) write $b=b\wedge (b\vee c)$ and use (ii)(a). (3) is proved 
similarly. For (4) first use (2) to get $\{a,b\}\prec c\vee d$, and then use
(ii)(a). (5) is proved similarly.
$\eop_{\ref{basic}}$
\end{Proof}

\begin{Lemma}\label{involution} Suppose that $B$ is a spil. Then $B$ satisfies:

{\noindent (1)} 
for all $x,y$ and $z$ we have $x \wedge y' \prec z$  iff  $x \prec z \vee y$, and

{\noindent (2)} for all $x$ and $y$,
$y\vee y'\prec x\implies 1\prec x$. 
\end{Lemma}

\begin{Proof} (1) Suppose that $x \wedge y' \prec z$, so by (iii)(b) we have $x\prec z\vee y''
=z\vee y$.

{\noindent (2)} Suppose that
$y\vee y'\prec x$. We have by the properness of the involution that
$y\vee y'=y''\vee y'$ which is by De Morgan laws equal to $(y'\wedge y)'$.
Hence $1\wedge (y'\wedge y)'=(y'\wedge y)'\prec x$.
By (iii)(a) we have $1\prec x\vee (y'\wedge y)$. By (ii)(a) there are $x^+$ and $y^+$
such that $x^+\prec x$ and $y^+\prec y'\wedge y$ such that $1\prec x^+\vee y^+$.
By (iv)(a) $y^+\prec y'\wedge y$ implies that $y^+\prec 0$ so $1\prec x^+\vee 0=x^+$
by Lemma \ref{basic}(4). By transitivity we get $1\prec x$.
$\eop_{\ref{involution}}$
\end{Proof}

We now proceed to associate to every spil a compact Hausdorff space, in a manner
similar to the classical Stone representation theorem.
The main difference is that filters are
defined in  connection with the $\prec$ relation rather than the Boolean-algebraic
order $\le$ and that there are no complements.

\begin{Definition}\label{deffilters} Suppose that $B$ is a spil.

{\noindent (1)} For $A\subseteq B$ we define $\uparrow A\deq\{x\in B:\,(\exists
a\in A)\,a\prec x\}$.

{\noindent (2)} A $\prec$-{\em filter} $F$ on $B$ is a non-empty
subset of $B$ 
which is closed under (finite) meets and satisfies $F=\uparrow F$.

{\noindent (3)} A $\prec$-filter $F$ on $B$ is 
called {\em prime} iff for every finite
$M\subseteq B$ with $\bigvee M\in F$ we have that $a\in F$ for some $a\in M$. 

{\noindent (4)} ${\bf spec}(B)$ is the set of all prime $\prec$-filters
with the topology generated by the sets
\[
O_x\deq \{F\in {\bf spec}(B):\,x\in F\}.
\]
\end{Definition}

Note that a prime $\prec$-filter is not necessarily an ultrafilter in the
sense of containing every set or its complement, as there is no
complement to speak of-- the involution does not necessarily satisfy
$x\wedge x'=0$ for all $x$.
That is why ${\bf spec}(B)$ is
not necessarily isomorphic to a subspace of $2^B$ and in fact it is not
necessarily zero-dimensional.
Some basic properties of prime filters are given by the following

\begin{Lemma}\label{filters} Let $B$ be a spil. Then:
\begin{description}
\item{(1)} if $F$ is a prime $\prec$-filter on $B$ then $0\notin F$, and
$1\in F$,
\item{(2)} if $a,b\in B$ then
$O_{a\wedge b}=O_a\cap O_b$ and $O_{a\vee b}=O_a\cup O_b$,
\item{(3)} if $F$ is a prime $\prec$-filter on $B$ and $a\in F$ then $a'\notin F$,
\item{(4)} if $F$ is a prime $\prec$-filter on $B$, $a,b\in F$ and
for some $x$ we have $x\prec a$ and $x'\prec b$, then $a\in F$ or $b\in F$,
\item{(5)} if $F\neq G$ are two prime $\prec$-filters on $B$, there is
$a$ such that $a\in F$ and $a'\in G$ or $a'\in F$ and $a\in G$.
\end{description}
\end{Lemma}

\begin{Proof} (1) If $0\in F$ then $\bigvee\emptyset\in F$
so $F\cap\emptyset\neq\emptyset$ by primeness, a contradiction. Since $\emptyset
\subseteq F$ we have $\bigwedge\emptyset\in F$ so $1\in F$.

{\noindent (2)} If $F$ is a $\prec$-filter containing
both $a,b$ then it also contains $a\wedge b$ by the closure under meets.
If $F$ is a $\prec$-filter containing $a\wedge b$ then by $F=\uparrow F$ we get
that for some $x\in F$ the relation $x\prec a\wedge b$ holds. Then $x\prec a$
and $x\prec b$ by the axioms of a spil, and hence $a,b\in F$. This shows the
first equality. For the second equality, if $F\in {\bf spec}(B)$
and $a\vee b \in F$ then by the primeness of 
$F$ we have $a,b\in F$; hence $O_{a\vee b} \subseteq O_a
\cup O_b$. If $F\in O_a$ then $a\in F=\uparrow F$, so for some $c\in F$ we have
$c\prec a$. By Lemma \ref{basic}(2) we have $c\prec a\vee b$ and hence
$a\vee b\in \uparrow F=F$. This shows $O_a\subseteq O_{a\vee b}$ and similarly
$O_b\subseteq O_{a\vee b}$.

{\noindent (3)} Suppose otherwise and let $a, a'\in F$, hence $a\wedge a'\in F=
\uparrow F$. By axiom (iv)(b) we have $a\wedge a'\prec 0$
so $0\in F$, contradicting (1).

{\noindent (4)} By Lemma \ref{basic}(4) we have $x\vee x'\prec a\vee b$.
By axiom (iv)(b) we have $1\prec x\vee x'$ then by (1)
above and the transitivity of $\prec$ we get $a\vee b\in F$, and hence
$a\in F$ or $b\in F$.

{\noindent (5)} Suppose $F\neq G$ and say $a\in F\setminus G$ (if there is no
such $a$, then there is $a\in G\setminus F$ and that case is handled by symmetry).
Since $a\in F=\uparrow F$ there is $b\in F$ with $b\prec a$, and for the same
reason there is $c\in F$ with $c\prec b$. By transitivity we have $c\prec a$.
By Lemma \ref{basic}(4) it follows that $c\prec a\vee b$, so by axiom (iii)(b)
of a spil we have $c\wedge b'\prec a$. On the other hand, by Lemma \ref{basic}(5)
we have $c\wedge b\prec a$. Putting these two conclusions together and
using Lemma \ref{basic}(4) we have $c\wedge (b\vee b')\prec a$. Using axiom 
(iii)(a) we have $b\vee b'\prec a\vee c'$ and then by axiom (iv)(b)
this implies $1\prec a\vee c'$. By (1) of this Lemma
we have $a\vee c'\in G$ so by the primeness of $G$ we have $a\in G$
or $c'\in G$. Since $a\notin G$ we have $c\in G$.
$\eop_{\ref{filters}}$
\end{Proof}

To prove Theorem \ref{forward} below we need to assure Hausdorffness and compactness
of the resulting space. The former will follow by Lemma \ref{filters} and for the
latter we shall need the following lemmata.

\begin{Lemma}\label{pomocni} Suppose that $B$ is a spil and $A\subseteq B$. Then:

{\noindent (1)} $\uparrow(\uparrow A)=\uparrow A$ and,

{\noindent (2)} if $A$ is closed under meets then so is $\uparrow A$.
\end{Lemma}

\begin{Proof} (1) If $c\in \uparrow A$ then there is $a\in A$ with $a\prec c$, so by 
axiom (i)(b) of spils there is some $b$ such that $a\prec b$ and $b\prec c$. Then $b\in
\uparrow A$, so $c\in \uparrow(\uparrow A)$.

If $c\in \uparrow(\uparrow A)$ then there is $b\in
\uparrow A$ such that $b\prec c$, hence $a\in A$ such that $a\prec b$ and $b\prec c$.
Since $\prec$ is transitive we have that $c\in \uparrow A$.

{\noindent (2)} Let $b, d\in \uparrow A$, hence there are $a,c\in A$ such that
$a\prec b$ and $c\prec d$. Then by Lemma \ref{basic}(5) we have $a\wedge b\prec c\wedge d$
and since $a\wedge b\in A$ we conclude $c\wedge d\in \uparrow A$.
$\eop_{\ref{pomocni}}$
\end{Proof}

\begin{Lemma}\label{Zorn} Suppose that $B$ is a spil and $A\subseteq B$ is closed
under meets and satisfies that for no $x\in A$ do we have $x\prec 0$.
Then there is a prime filter $F$ containing $A$ as a subset.
\end{Lemma}

\begin{Proof} Let $\FF$ be given by
\[
\FF=\{F\subseteq B:\,A\subseteq F, 0\notin F\mbox{ and }F \mbox{ is a filter}\}.
\]
By the choice of $A$ we have $0\notin \uparrow A$, so by Lemma \ref{pomocni}(2) we have
$\uparrow A$ is closed under meets. By Lemma \ref{pomocni}(1) we have 
$\uparrow(\uparrow A)=\uparrow A$, so $A\in \FF$. Consequently $\FF\neq\emptyset$. Now we
observe the following 

\begin{Claim}\label{1} If $F\in \FF$ then $\uparrow (F\cup \{1\})\in \FF$.
\end{Claim}

\begin{Proof of the Claim} By Lemma \ref{pomocni} it suffices to check that $\FF\cup
\{1\}$ is closed under meets and does not contain $0$, which follows by the choice of $F$.
$\eop_{\ref{1}}$
\end{Proof of the Claim}

It is easily seen that $\FF$ is closed under $\subseteq$-increasing unions so by Zorn's
lemma there is a maximal element $F$ of $\FF$. We claim that $F$ is prime. By Claim \ref{1}(1)
and maximality we have that $1\in F$. Now we shall show that for all $p\in B$ either $p$
or $p'$ are in $F$ (not both as then $0\in F$). So suppose that
$p\in B$ is such that $p,p'\notin F$. The family $X=\uparrow(F\cup\{p\wedge q:\,q\in F\})$ is clearly
a set satisfying $X=\uparrow X$ that is closed under meets and is proper a superset of $F$.
By maximality
of $F$ we have that for some $q\in F$ the relation $p\wedge q\prec 0$ holds.
Similarly we can find $r\in F$ such that $p'\wedge r\prec 0$ holds. Applying axiom (iii)(b)
of a spil we obtain that $q\prec p'$ and $r\prec p''$, so $p\wedge q\prec p'\wedge p''$
by Lemma \ref{basic}(5), and hence by axiom (iv) of a spil, $q\wedge r\prec 0$,
which is a contradiction with the choice of $F$.

Now suppose that $M\subseteq B$ is finite such that $m=\bigvee M\in F$ but no $p\in M$ is in
$F$. Hence for all $p\in M$ we have $p'\in M$ and so $\bigwedge\{p':\,p\in M\}=m'\in F$.
But then $m\wedge m'\in F$, which contradicts axiom (iv) and the fact that $0\notin F$.
We have shown that $F$ is as required.
$\eop_{\ref{Zorn}}$
\end{Proof}

\begin{Theorem}\label{forward} Let ${\bf spec}(B)$ be as defined in
Definition \ref{filters}. Then ${\bf spec}(B)$ is a compact Hausdorff space
with $\{O_x:\,x\in B\}$ a base.
\end{Theorem}

\begin{Proof} Clearly every element of ${\bf spec}(B)$ is contained in some
$O_a$. It follows by Lemma \ref{filters}(2) that the family 
$\{O_a:\,a\in B\}$ indeed forms a base for a topology on ${\bf spec}(B)$.
Now we show that the topology is Hausdorff.

Suppose that $F\neq G$ are prime $\prec$-filters and
let $a\in F\Delta G$. By Lemma \ref{filters}(5) there is $a$ such $a\in F$ and
$a'\in G$, or vice versa. Let us say that $a\in F$. Then $F\in O_{a}$ and $G\in
O_{a'}$ and by Lemma \ref{filters}(3), the sets $O_a$ and $O_{a'}$ are disjoint.

Finally we need to show that ${\bf spec}(B)$ is compact. So suppose that $\{O(p):\,p\in A\}$
covers ${\bf spec}(B)$ but no finite subfamily does. By Lemma \ref{filters}(2)
we may assume that $A$ is closed under finite joins. By the choice of $A$ for all finite $M\subseteq A$
there is $F\in {\bf spec}(B)$ with $\bigvee M\notin F$. Fix such an
$M$ and let $q=\bigvee M$.
If for some $p\in F$ we have that $p\wedge q'\prec 0$ then $p\prec 0\vee q''=q$, 
so $q\in F$ as $F$ is a filter,
a contradiction. So for no $p\in F$ do we have $p\wedge q'\prec 0$ and in 
particular we cannot have $q'\prec 0$ by Lemma \ref{basic}(3). This means that
the family $\{p':\,p\in A\}$ is closed under meets (as $A$ is closed under joins) and
none of its elements is $\prec 0$. By Lemma \ref{Zorn} there is a prime filter $F$
that contains this family as a subset. By the choice of $A$ there is $p\in A$ such
that $F\in O(p)$. But then $p, p'\in F$ which contradicts Lemma \ref{filters}(3).
$\eop_{\ref{forward}}$
\end{Proof}

The idea behind the direction from the space to a spil
in the representation theorem is that the pairs of the form $(O,K)$
where
$O$ is open
and $K\supseteq O$ compact will replace the clopen sets in the
Stone representation. The relation $\prec$ will be a replacement 
for $\subseteq$ (so $\le$ in the Ba representation), so we shall have
$(O_0,K_0)\prec (O_1,K_1)$ iff $K_0\subseteq O_1$.

\begin{Theorem}\label{backward} Suppose that $X$ is a compact Hausdorff space.
Let $\TT=\TT_X$ denote the set of all open subspaces of $X$ and $\KK=\KK_X$
the set of all compact subspaces of $X$. We define
\begin{itemize}
\item
$B\deq\{(O,K):\,O\in\TT, K\in \KK, O\subseteq K\}$,
\item
$(O_0,K_0)\vee (O_1,K_1)\deq (O\cup O_1, K_0\cup K_1)$,
\item
$(O_0,K_0)\wedge (O_1,K_1)\deq (O_0\cap O_1, K_0\cap K_1)$,
\item
$0\deq (\emptyset,\emptyset), 1\deq (X,X)$,
\item
$(O_0,K_0)\prec (O_1,K_1)\iff K_0\subseteq O_1$,
\item
$(O,K)'\deq (X\setminus K, X\setminus O)$.
\end{itemize}
Then $\langle B, \vee, \wedge, 0,1,\elementary,'\rangle$ is a spil such that
${\bf spec}(B)$ is homeomorphic to $X$.
\end{Theorem}

\begin{Proof} It is clear that $\langle  B, \vee, \wedge, 0,1\rangle$
is a distributive bounded lattice,
as well as that $\prec$ is transitive. Since $X$ is
compact Hausdorff it is normal so the operation $\prec$ is indeed interpolating.
The second axiom from the list in Definition \ref{sba}
is easily seen to hold by the definition of $\wedge$ and $\vee$. Let us
consider axiom (iii). 

The involution is clearly proper. For part (b)
suppose that $(O_0,K_0)\wedge (O_1, K_1)\prec
(O_2, K_2)$, so $K_0\cap K_1\subseteq O_2$. We have $(O_1, K_1)'=
(X\setminus K_1, X \setminus O_1)$ so $(O_2, K_2)\vee (O_1, K_1)'=
(O_2\cup (X\setminus K_1), K_2\cup (X\setminus O_1))$. Since $K_0
\subseteq O_2\cup (X\setminus K_1)$ we obtain that $(O_0,K_0)\prec 
(O_2, K_2)\vee (O_1, K_1)'$, as required. The remaining direction of
the axiom is proved similarly. De Morgan laws clearly hold.

For axiom (iv), if $(O,K)\prec (U,H)\wedge (X\setminus H, X \setminus U)$ then
since $U\subseteq H$ we have $X\setminus U\supseteq X\setminus H$ and hence
$U\cap (X\setminus H)=\emptyset$ (as a side note observe that it does not
necessarily follow
that $H\cap (X\setminus U)=\emptyset$). Since $(O,K)\prec (U,H)$ we have $K=
\emptyset$, so $O=\emptyset$ and clearly $(O,K)\prec (\emptyset,\emptyset)$. 

This shows that $B$ is a spil and we have to verify that
$X$ is homeomorphic to ${\bf spec}(B)$. To this end let us define for $x\in X$
the set $F_x=\{(O,K)\in B:\,x\in O\}$.

\begin{Claim}\label{defined} Each $F_x$ is an element of ${\bf spec}(B)$.
\end{Claim}

\begin{Proof of the Claim} Let $x\in X$. Since $(X,X)\in F_x$ we have that
$F_x\neq\emptyset$. It is clear that $F_x$ is closed under meets, so $F_x$
is a filter. Suppose that $(\bigcup_{i<n}O_i, \bigcup_{i<n}K_i)\in F_x$,
where each $(O_i,K_i)\in B$. Hence $x\in \bigcup_{i<n}O_i$ so there is
some $i<n$ such that $x\in O_i$ and so $(O_i,K_i)\in F_x$.
$\eop_{\ref{defined}}$
\end{Proof of the Claim}

Let $g$ be the function associating $F_x$ to $x$. We claim that $g$ is a
homeomorphism between $X$ and ${\bf spec}(B)$. If $x\neq y$ then there is
$O$ open containing $x$ and not containing $y$. Hence $(O,X)\in F_x\setminus F_y$
and hence $F_x\neq F_y$. So $g$ is 1-1.

Suppose that $F\in {\bf spec}(B)$ and let $\KK=\{K:\,(\exists O)(O,K)\in F\}$.
Since this is a centred family of compact sets its intersection is non-empty,
so let $x\in\bigcap\KK$. We claim that $F=F_x$. If not, then there is
$a=(O,K)\in F_x$ such that $a'\in F_x$ (by Lemma \ref{filters}(5) and the fact that
the involution is proper in $B$). But then $x\in O$ and hence $x\notin
X\setminus O$, contradicting the assumption that $a'=(X\setminus K, X\setminus O)\in
F$. Hence $g$ is bijective.

Suppose that $U$ is basic open in ${\bf spec}(B)$ so $U=O_a$ for some $a=(O,K)$.
Then 
\[
g^{-1}(O_a)=\{x:\,F_x\in O_a)\}=\{x:\,a\in F_x\}=\{x:\,x\in O\}=O,
\]
so open in $X$. Hence $g$ is continuous.

Finally, if $O$ is open in $X$ then $g``O=\{g(x):\,x\in O\}=\{F_x:\,x\in O\}$.
If $U$ is open $\subseteq O$ and $K$ is a compact superset of $U$ then
if $F=O_{(U,K)}$, $F=F_x$ for some $x\in U$, as follows
from the argument showing surjectivity of $g$.
Hence $O_{(U,K)}\subseteq\{F_x:\,x\in O\}$, which
shows that $=\{F_x:\,x\in O\}$
contains $\bigcup\{O_{(U,K)}:\,U\mbox{ open}\subseteq O, K\mbox{ compact}\supseteq
U\}$. In fact we claim that these two sets are equal, which shows that $g$ is an
open mapping and hence a homeomorphism. So let $x\in O$ and $(U,K)\in F_x$.
Hence $(O\cap U,K)\in F_x$ and so $F_x\in O_{(U,K)}$.
$\eop_{\ref{backward}}$
\end{Proof}

We finish this section by explaining the use of the word ``strong" in the name for a spil.
In the terminology of \cite{JS}, proximity lattices are structures that satisfy 
the axioms of a spil but without the involution, and such structures are
called strong if they in addition satisfy the following axioms
\begin{description}
\item{(A)}\label{iv}
for all $a,x,y\in B$
\[
x\wedge y\prec a\implies (\exists x^+,y^+\in B)\, x\prec x^+, y\prec y^+
\,\,\&\,\,
x^+\wedge y^+\prec a;
\]
\item{(B)}\label{v}
for all $a,x,y\in B$
\[
a\prec x\vee y\implies (\exists x^+,y^+\in B)\, x^+\prec x, y^+
\prec y \,\,\&\,\,
a\prec x^+\vee y^+;
\]
\end{description}
Note 
that $\prec$ is not necessarily reflexive in a spil hence axioms (A) and (B) are not
trivially met. We shall however demonstrate that every spil
satisfies them.
 
\begin{Claim}\label{strong}
Suppose $B$ is a spil.
Then axioms (A) and (B) above are satisfied.
\end{Claim}

\begin{Proof of the Claim} Let us first show (A), so suppose that 
$x\wedge y\prec a$.  Then by the interpolating property of $\prec$ there
is $b$ such that $x\wedge y\prec b\prec a$. By axiom (iii)(b) of a spil
this gives $x\prec b\vee y'$. Similarly we obtain $y\prec b\vee x'$.
Letting $x^+=b\vee y'$ and $y^+= b\vee x'$ we have $x^+\wedge y^+=
b\wedge (x'\vee y')$. Since $b\prec a$, by Lemma \ref{basic}(3) we
have $b\wedge (x'\vee y')\prec a$, hence $x^+$ and $y^+$ are as required.

(B) is shown similarly.
$\eop_{\ref{strong}}$
\end{Proof of the Claim}

\section{Lattices}\label{lattices} Here we show that spils occur naturally in the
context of sublattices of Boolean algebras.

\begin{Definition}\label{latproduct} Suppose that $\mathfrak A$ is a fixed Boolean algebra
with Boolean operations $+,\cdot$ and $-$, and the induced order $\le$,
and that $L$ is a sublattice of $\mathfrak A$ satisfying
the following {\em normality} condition: 
for all $a, b\in L$ satisfying $a\cdot b=0_{\mathfrak A}$ there are
$c,d\in L^c$ satisfying $a\le c, b\le d$ and $c\cdot d=0_{\mathfrak A}$.

We define {\em the spil induced by $L$} by letting
\[
B=\{(u,k):\,u\in L^c, k\in L\mbox{ and }u\le k\}
\]
and endow it with the following operations:
\begin{itemize}
\item $(u,k)\wedge (v,h)=(u\cdot v, k\cdot h)$,
\item $(u,k)\vee (v,h)=(u + v, k+ h)$,
\item $(u,k)'=(-k,-u)$
\end{itemize}
and
the relation $(u,k)\le (v,h)$ iff  $k\le v$. We let
$1=(1_{\mathfrak A}, 1_{\mathfrak A})$ and $0=(0_{\mathfrak A}, 0_{\mathfrak A})$.
\end{Definition}

\begin{Theorem}\label{induced}
Suppose that $L$ is as in Definition \ref{latproduct}. Then

{\noindent (1)} the spil $B$ induced by $L$ is a spil and

{\noindent (2)} the lattice of closed subsets of ${\bf spec}(B)$ is isomorphic to
$L$.
\end{Theorem}

\begin{Proof} (1) Clearly $B$ is a bounded distributive
lattice with the 0 and 1 as specified. We check the rest of the axioms of Definition
\ref{sba}. 

Clearly $\prec$ is transitive. Checking
that the relation $\prec$ is interpolating uses the normality
of $L$. Suppose that 
$(u,k)\prec (v,h)$ holds, hence $k\le v$ and hence $k, -v$ are disjoint elements of $L$.
Let $w\ge k$ and $z\ge -v$ be disjoint elements of $L^c$. Then $(w, -z)\in B$
satisfies $(u,k)\prec (w,-z)\prec (v,h)$.

Axiom (ii) of a spil follows by the corresponding properties of the Boolean algebra
$\mathfrak A$. Similarly for axioms (iii) and (iv).

{\noindent (2)} Let $B^\ast$ be the spil consisting of pairs $(O,K)$ of pairs of
open and compact subsets of ${\bf spec}(B)$ such that $O\subseteq K$. By the
representation theorem in \S\ref{spils} we have that ${\bf spec}(B^\ast)$ and ${\bf spec}(B)$
are homeomorphic and this induces an isomorphism between $B$ and $B^\ast$ given by
the identity function. Hence $L$ is isomorphic to the lattice of closed, equivalently,
compact, subsets of ${\bf spec}(B)$.
$\eop_{\ref{induced}}$
\end{Proof}

Theorem \ref{induced} and the representation theorem from \S1 give a representation theorem
for normal lattices due to Plebanek in \cite{Grz} where the theorem is proved directly.
He showed there also that the lattice being in addition connected (so $L\cap L^c=\{0,1\}$)
and interpolating (for all $a\in L\setminus\{0,1\}$ there are $l,u$ such that
$a\cdot l=0_{\mathfrak A}$ and $a+ u=1_{\mathfrak A}$) imply that ${\bf spec}(B)$ defined
as above is connected.

Which theorem to use in constructions of course depends on the context: when one is working
inside of a fixed Boolean algebra then one might prefer to construct a lattice, while
if no ambient Boolean algebra is specified then a spil might be easier to construct.

\end{document}